\documentclass[12pt]{article}
\usepackage{amsmath,amssymb}
\usepackage{epsf}
\numberwithin{equation}{section}

\begin{document}
\newcommand{\derl}{\partial_L}
\newcommand{\derr}{\partial_R}
\newtheorem{prop}{Proposition}

\def\tuc{\bf}
\def\tr{{\rm tr}\,}
\def\ba{\begin{array}}
\def\ea{\end{array}}
\def\be{\begin{equation}}
\def\ee{\end{equation}}
\def\bea{\begin{eqnarray}}
\def\eea{\end{eqnarray}}
\def\wc{\bar{w}}
\def\ld{\ldots}
\def\ie{i.e.}
\def\C{{\mathbb C}}
\def\Z{{\mathbb Z}}
\def\R{{\mathbb R}}
\def\cf{{\mathcal{F}}}
\def\cs{{\mathcal{S}}}
\def\cl{{\mathcal{L}}}
\def\cx{{\mathcal{X}}}
\def\cz{{\mathcal{Z}}}
\def\su{{\mathfrak s \mathfrak u}}
\def\fd{{f^\dagger}}
\def\dg{{^\dagger}}
\def\fdo{{f^{0\dagger}}}
\def\1{{\bf 1}}
\def\nn{\nonumber}

\title{Surfaces associated with sigma models and the sine--Gordon equation}

\author{
A.~M. Grundland\thanks{email address: grundlan@crm.umontreal.ca}
\\
Centre de Recherches Math{\'e}matiques, Universit{\'e} de Montr{\'e}al, \\
C. P. 6128, Succ.\ Centre-ville, Montr{\'e}al, (QC) H3C 3J7, Canada \\
Universit\'{e} du Qu\'{e}bec, Trois-Rivi\`{e}res CP500 (QC) G9A 5H7, Canada \\[2mm]
{\rm and} \\[2mm]
L. \v{S}nobl\thanks{email address: Libor.Snobl@fjfi.cvut.cz}
\\
Centre de Recherches Math{\'e}matiques, Universit{\'e} de Montr{\'e}al, \\
C. P. 6128, Succ.\ Centre-ville, Montr{\'e}al, (QC) H3C 3J7, Canada \\
and \\
Faculty of Nuclear Sciences and Physical Engineering, \\
Czech Technical University, \\
B\v rehov\'a 7, 115 19 Prague 1, Czech Republic\\} \date{}

\maketitle

\abstract{We present a unified method of construction of surfaces associated with Grassmannian sigma models, expressed in terms of an orthogonal projector. This description leads to compact formulae for structural equations of two-dimensional surfaces immersed in the $su(N)$ algebra. In the special case of the $\C P^1$ sigma model we obtain constant negative Gaussian curvature surfaces. As a consequence, this leads us to an explicit relation between the 
$\C P^1$ sigma model and the sine--Gordon equation.} 
\smallskip

\smallskip


\section{Introduction}

Description of the behaviour of surfaces immersed in $\R^n$ in connection with integrable systems lead in most cases
to fundamental forms where the coefficients satisfy the Gauss--Weingarten and Gauss--Codazzi--Ricci equations. The study of general properties of these equations and the methods of solving them is a rapidly developing area of modern mathematics. A broad review of recent developments in this subject can be found e.g. in \cite{Hel1,Hel2} (and references therein).

The idea of inducing surfaces in three--dimensional Euclidean space from the solutions of two--dimensional linear problems is a very old one. It originates from the work of K. Weierstrass \cite{Wei} and A. Enneper \cite{Enn} one and half centuries ago. This subject has been extensively researched (see e.g. \cite{Ken1,Ken2,Hof1,Hof2,Oss}).  Until very recently, the immersion of two--dimensional surfaces obtained through the Weierstrass formula had only been known in low dimensional Euclidean spaces \cite{Bud,Sas}. The relationship between the immersion of surfaces and Dirac--type systems has been developed over the last decades by many authors \cite{Kon,Kon1,Fri} with a purpose of extension of the applicability of the immersion formula and consequently of constructing more diverse type of surfaces than those obtained previously by the classical approach. For instance surfaces associated with sigma models provide us with a rich class of geometric objects. The rich character of this formulation makes immersion formula an interesting object of study and various special types of surfaces have been investigated (\cite{Pin}). 

This paper is a follow--up of the results obtained in \cite{Grusno,Grusnosym,GrusnoGr} and is concerned with two--dimensional surfaces immersed in multi--dimensional Euclidean spaces obtained from solutions of Grassmannian sigma models defined on Minkowski space. The heart of the matter is that the equations defining the immersion are formulated directly in terms of matrices taking their values in the Lie algebra $su(N)$. The main advantages of this procedure is that the group analysis of the immersion makes it possible to construct regular algorithms for finding certain classes of surfaces without referring to any additional considerations. The proposed method proceeds directly from the given sigma model.

The task of finding the surfaces is facilitated by the group properties of these models. Our main goal is to provide a self--contained comprehensive approach to this subject. For this purpose we formulate the structural equations for the immersion expressed in the Cartan's language of moving farmes, fundamental forms, the Gauss curvature and the mean curvature vector, using an orthogonal projector satisfying the Euler--Lagrange equations of the given sigma model. The main advantage of the projector approach is that there are no gauge degrees of freedom in this description of the Grassmannian sigma models, compared to the use of equivalence classes. Such a description leads to much simpler formulae and allows us to write in closed form quantities which were previously too complicated to be presented. 

Secondly, we return to study the simplest case of this construction, namely the $\C P^{1}$ sigma model. In this particular case, the resulting surfaces have negative constant Gaussian curvature. Therefore, we may construct and study the corresponding solutions of the sine--Gordon equation. It turns out that the relation doesn't necessarily involve the construction of surfaces, i.e. there is a direct reduction from the $\C P^{1}$ sigma model to the the sine--Gordon equation, provided certain regularity conditions on the $\C P^1$ solution are met.

In order to illustrate this relation we consider solutions of the $\C P^{1}$ sigma model obtained in different ways -- a symmetry reduction, B\"acklund transformation -- and construct both the corresponding surfaces and the solutions of the sine--Gordon equations. 

\section{Grassmannian sigma models and their equations of motion}\label{SecGrsm}

We generalize and simplify our treatment of sigma models on complex Grassmannian manifolds defined on Minkowski space. 

The starting point of our generalization lies in the realization that most of the properties of associated surfaces 
were described using a projector $P$,
$$ P^2=P, \; P^\dg =P,$$
which satisfies the equations of motion in the form
\be\label{ELP}
[\derl \derr P, P] =0.
\ee
Therefore, the description of the Grassmannian sigma models in terms of the projector 
on the corresponding subspace defining the element of the Grassmannian manifold
\be\label{GrP} G(m,n) = \{ {\rm m-dimensional \; subspaces \; of \; } \C^{N} \} \ee
is more natural in our context than the description using equivalence classes of the elements of 
\be\label{Gmn}
G(m,n) = \frac{SU(N)}{S(U(m) \times U(n))}, \; \; N=m+n.
\ee

In this formalism, the solution of the model is described 
as a map 
$$ P:\Omega \rightarrow {\rm Aut}(\C^{N}), \qquad P^\dg=P, \qquad P^2=P$$
where $\Omega \subset \R^2$ will be assumed to be a connected and simply 
connected domain of definition of the model.
The solution is required to be a stationary point of the action  
\be\label{GRmodel}
{\cs} = \int_\Omega \tr \{ \derl P . \derr P \} {\rm d}\xi_L {\rm d}\xi_R
\ee
where $\xi_L,\xi_R$ are the light-cone coordinates on the Minkowski space $R^2$, i.e. the metric on $R^2$
is written as
\be\label{Minkmetr}
 {\rm d}s^2= {\rm d}\xi_L {\rm d}\xi_R. 
\ee
We shall denote by $\derl$ and $\derr$ the derivatives with respect to $\xi_L$ and $\xi_R$, respectively.

Because the hermitian matrix $P$ is subject to the constraint $$P^2=P$$
we have to introduce the Lagrange multiplier $\lambda=\lambda^{\dg} \in {\rm Aut}(\C^{N})$ into the action (\ref{GRmodel})
\be\label{LagrMult}
{\cs} = \int_\Omega \tr \{ \derl P . \derr P + \lambda . (P^2-P) \} {\rm d}\xi_L {\rm d}\xi_R.
\ee
By the variation of the action (\ref{LagrMult}) we get
\bea
\nn \delta \lambda: & & \; \; \; P^2-P  = 0, \\
\label{varP} \delta P : & & \; \; \; 2 \derl \derr P + \lambda.P +P.\lambda-\lambda =0.
\eea
In order to eliminate the Lagrange multiplier $\lambda$ from the Euler--Lagrange equations we 
multiply (\ref{varP}) from the left (right) by $P$ and subtract them, obtaining (\ref{ELP})
$$ [\derl \derr P,P] =0 $$
For later reference let us note that an equivalent version of (\ref{ELP}) reads
\be\label{ELPeqcons}
\derl [\derr P, P] + \derr [\derl P,P] =0
\ee
and that as a differential consequence of $P^2=P$ one has
\be\label{Pprops}
 \partial_{D} P = \partial_{D} P P + P \partial_{D} P, \; P \partial_{D} P P =0, \; D=L,R . 
\ee
Let us also introduce a notation for the trace of a product of two derivatives of $P$
\be
p_{B_1\ldots B_k|D_1\ldots D_l} = \tr \left( \partial_{B_1\ldots B_k} P. \partial_{D_1\ldots D_l} P \right),
\ee
where $k,l>0, \; B_1,\ldots,B_k,D_1,\ldots,D_l = L,R.$ 
Note that $p_{B_1\ldots B_k|D_1\ldots D_l}=p_{D_1\ldots D_l|B_1\ldots B_k}$.

Finally, let us mention that projectors satisfying (\ref{ELP}) may also arise in other models 
or applications which may not be naturally interpreted in terms of Grassmannian manifolds. 
Nevertheless, the construction of surfaces associated with them can be performed in the same way.

\section{Surfaces obtained from Grassmannian sigma model}\label{Secsurf}

Let us now discuss the analytical description of a two--dimensional surface $\cf$ 
immersed in the $\su(N)$ algebra, associated with the projector (\ref{ELP}).
Firstly, we shall construct an exact $\su(N)$--valued 1--form whose ``potential'' 0--form defines the surface $\cf$. 
Next, we shall investigate the geometric characteristics of the surface $\cf$. 

Let us introduce the scalar product
$$ ( A , B ) = -\frac{1}{2} \tr A.B $$
on $\su(N)$ and identify the $(N^2-1)$--dimensional Euclidean space with the $\su(N)$ algebra
\be\label{ident} \R^{N^2-1} \simeq \su(N). 
\ee
We denote 
\be\label{mlmr} M_L = [\derl P,P], \qquad M_R = [\derr P,P]. \ee
It follows from (\ref{ELPeqcons}) that 
\be\label{ELeqcons2}
\derl M_R+ \derr M_L =0 .
\ee
We identify tangent vectors to the surface $\cf$ as follows
\be\label{dlxdrx}
 {\cx}_L=M_L, \; \; {\cx}_R=-M_R .
\ee
Equation (\ref{ELeqcons2}) implies the existence of a closed $\su(N)$--valued 1--form on $\Omega$
$$ {\cx} = {\cx}_L {\rm d} \xi_L + {\cx}_R {\rm d} \xi_R, \; \;  {\rm d} {\cx} =0 .$$
Because ${\cx}$ is closed and $\Omega$ is by assumption connected and simply connected, ${\cx}$ is also exact. In other words,
there exists a  well--defined $\su(N)$--valued function $X$ on $\Omega$  that ${\cx} = {\rm d} X$. 
The matrix function $X$ is unique up to addition of any constant element of $\su(N)$ and 
we identify the components of $X$ with the coordinates of the sought--after surface $\cf$ in $\R^{N^2-1}$. 
Consequently, we get 
\be\label{cldx}
\derl X= {\cx}_L,\ \derr X= {\cx}_R.
\ee
The map $X$ is called the Weierstrass formula for immersion. 
In practise, the surface $\cf$ is found by integration
\be\label{surfaceN}
\cf: \;  X(\xi_L,\xi_R) = \int_{\gamma(\xi_L,\xi_R)} {\cx}
\ee
along any curve $\gamma(\xi_L,\xi_R)$ in $\Omega$ connecting the point $(\xi_L,\xi_R)\in \Omega$ with an arbitrary chosen point $(\xi_L^0,\xi_R^0)\in \Omega$.

We should investigate the behaviour of the constructed surface under the known symmetries of 
(\ref{ELP}). The equation (\ref{ELP}) is invariant under the following change of independent 
variables (i.e. the conformal transformation)
\be\label{conf}
 \xi_L \longrightarrow f(\xi_L), \; \; \xi_R \longrightarrow g(\xi_R).
\ee
Since the surface $\cf$ is written in terms of an integral of a one--form, such a transformation amounts only to a reparametrization of the surface; as a geometric object it remains the same. Another symmetry is the transformation 
\be
P \longrightarrow U P U^\dg, \; \; U \in U(N).
\ee
The only effect of such a transformation on the surface is a rotation in $\R^{N^2-1}$, so again the geometry of the surface is unchanged. Therefore, the surface $\cf$ associated with a solution of (\ref{ELP}) characterizes the symmetry equivalence class of solutions of (\ref{ELP}).

By computation of the scalar products of ${\cx}_{B}.{\cx}_{D}, \ B,D = L,R$ we find the components of the induced metric on the surface 
${\cf}$
\bea\label{metric} 
 G   & = & \left( \begin{array}{cc} J_{L}, & G_{LR} \\ G_{LR}, & J_{R} \end{array}  \right) = \\ 
 \nonumber & = & \frac{1}{2} \left( \begin{array}{cc} p_{L|L}  & - p_{L|R}  \\ 
- p_{L|R} & p_{R|R} \end{array}  \right). 
\eea
As a consequence of (\ref{ELP}) we find
$$ \derr J_{L} = p_{LR|L} = 2 \tr \left( \partial_{LR} P P \derl P \right)   = 2 \tr \left( \partial_{LR} P P \derl P P \right) =0$$
using the cyclic property of the trace and along with (\ref{Pprops}). Similarly
$$ \derl J_{R} =0.$$
The first fundamental form of the surface $\cf$ takes the compact form
\bea
\nonumber I   &= &   J_L ({\rm d}\xi_L)^2 -  2 G_{LR} {\rm d}\xi_L {\rm d}\xi_R +  J_R ({\rm d} \xi_R)^2 \\
& = & (\delta_{B,D}-\frac{1}{2}) p_{B|D} {\rm d}\xi_B {\rm d}\xi_D
\label{1stff}
\eea
where summation over repeated indices $B,D=L,R$ applies and $\delta_{B,D}=1$ if $B=D$ and 0 otherwise.
It can be shown using the Schwarz inequality that such a first fundamental form $I$ is positive, and then investigated under which conditions it is positive definite \cite{GrusnoGr}, i.e. when the surface is, at least locally, well defined.

It is useful to note that conformal transformations of independent variables change the metric (\ref{Minkmetr}) on $\R^2$ but leave invariant the Euler--Lagrange equations (\ref{ELP}. Since the metric (\ref{Minkmetr}) is no longer needed in the following we feel free to use such a transformation (\ref{conf}) to bring the solution of (\ref{ELP}) to an equivalent solution (outside singular points where $J_L . J_R=0$ and consequently the tangent vectors $\derl X,\derr X$ are linearly dependent) such that 
\be\label{normalization} J_L = \frac{1}{2} p_{L|L} =1, \; \; J_R = \frac{1}{2} p_{R|R} =1 .\ee
Such transformation is for a given solution of (\ref{ELP}) expressed in terms of quadratures. The first fundamental form now takes a particularly simple form 
\bea\label{Iinconf}
I   &= &   ({\rm d}\xi_L)^2 -  p_{L|R} {\rm d}\xi_L {\rm d}\xi_R +  ({\rm d} \xi_R)^2,
\eea
i.e. the surface is described in Chebyshev coordinates. In the following we shall always assume that we chose our independent coordinates in this way. This assumption allows us to write a lot of expressions in a closed form; otherwise they would be too complicated to be presented here.  

Using (\ref{metric}) we can write the formula for the scalar curvature \cite{DoC} as
\bea\label{Gaussian}
K  & = &  2 \left(  \frac{p_{LR|LR}-p_{LL|RR}}{4- (p_{L|R})^2} - \frac{p_{LL|R} \ p_{L|RR}  \ p_{L|R}}{(4 - (p_{L|R})^2 )^2}\right).
\eea

\section{The moving frames and the Gauss--Weingarten equations}\label{SecGW}

Now we may formally determine a moving frame on the surface $\cf$ and write the Gauss--Weingarten equations. Let $P$ be a solution of (\ref{ELP}) such that ${\rm det}(G)$ is not zero in a neighbourhood of a regular point $(\xi_L^0,\xi_R^0)$ in $\Omega$, so that we can assume (\ref{normalization}). Assume also that the surface $\cf$ (\ref{surfaceN}) associated with these equations is described by the moving frame on $\cf$
$$\vec \tau=({\cx}_L, {\cx}_R, n_{3}, \ldots, n_{N^2-1})^T,$$
where the vectors ${\cx}_L, {\cx}_R, n_{3}, \ldots, n_{N^2-1}$ are identified with matrices as in (\ref{ident}) and satisfy the normalization conditions
$$ ({\cx}_L,{\cx}_L)=1, \, ({\cx}_L,{\cx}_R)= -\frac{1}{2} p_{LR}, \, ({\cx}_R,{\cx}_R)=1, $$
\be\label{normnorm} ({\cx}_L, n_{k})=({\cx}_R, n_{k})=0, \, (n_j, n_{k})=\delta_{jk}. \ee
The moving frame satisfies the Gauss--Weingarten equations
\begin{eqnarray}
\nonumber \derl {\cx}_L & = & A^L_L {\cx}_L + A^L_R {\cx}_R + Q^L_j n_j, \\
\nonumber \derl {\cx}_R & = & H_j n_j, \\ 
\nonumber \derl n_j & = & \alpha^L_j {\cx}_L + \beta^L_j {\cx}_R +s^L_{jk} n_k, \\
\nonumber \derr {\cx}_L & = & H_j n_j, \\ 
\nonumber \derr {\cx}_R & = & A^R_L {\cx}_L + A^R_R {\cx}_R + Q^R_j n_j, \\
\derr n_j & = & \alpha^R_j {\cx}_L + \beta^R_j {\cx}_R +s^R_{jk} n_k \label{gweqN},
\end{eqnarray}
where $s^L_{jk}+s^L_{kj}=0$, $s^R_{jk}+s^R_{kj}=0,$ $j,k=3,\ldots, N^2-1$,
$$ \alpha^L_j = - 2\frac{ p_{L|R} H_j  + 2 Q^L_j }{  4- (p_{L|R})^2 }, \, \, \, 
\beta^L_j =  - 2 \frac{   p_{L|R} Q^L_j  +  2 H_j }{ 4- (p_{L|R})^2 },$$ 
$$ \alpha^R_j = - 2 \frac{ p_{L|R}  Q^R_j +  2 H_j}{ 4- (p_{L|R})^2 }, \, \,  \, 
\beta^R_j = - 2\frac{  p_{L|R}  H_j + 2Q^R_j }{ 4- (p_{L|R})^2 },$$ 
and $A^L_L,A^L_R,A^R_L,A^R_R$ are the Christoffel symbols of the second kind
given by
\begin{eqnarray}
\nonumber A^L_L  & = & \Gamma^{L}_{LL} = \frac{ - p_{L|R} \ p_{LL|R}}{ 4- (p_{L|R})^2 } ,\\
\nonumber A^L_R  & = & \Gamma^{L}_{RR} = \frac{ - 2 p_{LL|R}}{ 4- (p_{L|R})^2 } ,\\
\nonumber A^R_L  & = & \Gamma^{R}_{LL} = \frac{ - 2 p_{RR|L}}{ 4- (p_{L|R})^2 } ,\\
A^R_R  & = & \Gamma^{R}_{RR} = \frac{ - p_{L|R} \ p_{RR|L}}{ 4- (p_{L|R})^2 }. 
\end{eqnarray}

%

The explicit form of the coefficients $H_j,Q^D_j$ (where $D=L,R$; $j=3,\ldots, N^2-1$) depends on 
the chosen orthonormal basis $\{ n_3, \ldots, n_{N^2-1} \}$ of the space normal to the surface $\cf$ 
at the point $X(\xi_L^0,\xi_R^0)$. They are not completely arbitrary, since using (\ref{ELP}) and (\ref{Pprops}) we find that they are restricted by the condition
$$ \left( \derl {\cx}_L, \derl {\cx}_R \right)=\left( \derr {\cx}_R, \derl {\cx}_R \right)=0. $$

The derivation of the Gauss--Weingarten equations is almost the same as in \cite{GrusnoGr}, only in better notation, therefore we will not present it here and refer the interested reader to \cite{GrusnoGr}. Due to the current notation using $P$ and Chebyshev coordinates the formulae presented here are significantly simpler.

An example of a moving frame of the surface $\cf$ can be constructed as follows.
Let $P$ be a solution of (\ref{ELP}). Taking into account that
\be\label{traceinv} \tr(A)= \tr(\Phi A \Phi^{\dagger}), \; \;  \Phi \in SU(N), \ee
we may employ the adjoint representation of the group $SU(N)$ in order to bring $ \derl X,\derr X, n_a$ to its simplest possible form. We shall request that $\Phi(\xi_L,\xi_R)$ diagonalizes $P(\xi_L,\xi_R)$, i.e.
\be\label{Phireq} P(\xi_L,\xi_R)  =  \Phi(\xi_L,\xi_R) {\rm diag} (0,\ldots,0,1,\ldots,1) \Phi^\dg(\xi_L,\xi_R). \ee

The transformed derivatives of $X$ have the following block matrix structure (for the sake of brevity we suppress the dependence of $\Phi,P$ etc. on $\xi_L,\xi_R$).
\begin{eqnarray}\label{dz0}
\nonumber  \partial_D^\Phi X & \equiv & \Phi^\dagger \partial_D X \Phi  = \left( \ba{cc} 0 & \partial_{D}^\Phi P \\ -(\partial_{D}^\Phi P)^\dg & 0 \ea \right).
\end{eqnarray}
where $\partial_{D}^\Phi P$ is defined by the equality. 

Let us choose an orthonormal basis in $\su(N)$ in the following form
\begin{eqnarray}
  (A_{jk})_{ab} & = &  i (\delta_{ja} \delta_{kb} + \delta_{jb} \delta_{ka} ), \; \; 1\leq j<k\leq N, \\
\nonumber  (B_{jk})_{ab} & = &   (\delta_{ja} \delta_{kb} - \delta_{jb} \delta_{ka} ), \; \; 1\leq j<k\leq N, \\
\nonumber (C_{p})_{ab} & = & i \sqrt{\frac{2}{p(p+1)}} \left( \sum_{d=1}^p \delta_{da} \delta_{db} - p \delta_{p+1,a} \delta_{p+1,b} 
\right), \; \;  1 \leq p \leq N-1.
\end{eqnarray}
When a solution $\Phi$ of (\ref{Phireq}) is known, the construction of the moving frame proceeds as follows.
One finds, using the Gramm-Schmidt orthogonalization procedure, the orthonormal vectors 
$$ \tilde A_{aj}, \tilde B_{aj}, \; a \leq m, \; j > m, \; a+j > m+2 $$
satisfying the conditions
$$ (\partial_{D}^\Phi X, \tilde A_{aj}) =0, \, (\partial_{D}^\Phi X, \tilde B_{bj}) =0 $$
and 
\be
{\rm span} (\partial_{D}^\Phi X,\tilde A_{aj}, 
\tilde B_{aj})_{D=L,R, \ a \leq m}^{j > m, \; a+j > m+2}  = 
{\rm span} (A_{aj},B_{aj})_{a \leq m\; < j},
\label{GrammSch}
\ee
and then identifies the remaining tilded and untilded matrices
$$ \tilde A_{jk}= A_{jk}, \ \tilde B_{jk}= B_{jk}, \ \tilde C_{p} = C_{p}, $$
where $a,j \leq m$ or $a,j > m$ and $p \leq N-1.$
Consequently, $\partial_{D}^\Phi X, \tilde A_{jk}, \ \tilde B_{jk}, \ \tilde C_{p}$
satisfy the normalization conditions like (\ref{normnorm}). By the invariance of trace under unitary transformation
 (\ref{traceinv}) we find that the moving frame of the surface $\cf$ in the neighbourhood $\Upsilon$ of point $X^0=X(\xi_L^0,\xi_R^0)$ 
\begin{eqnarray}
\nonumber \derl X(\xi_L,\xi_R)  & = & \Phi(\xi_L,\xi_R) \derl^\Phi X(\xi_L,\xi_R) \Phi^\dagger(\xi_L,\xi_R), \\ 
\nonumber \derr X(\xi_L,\xi_R)  & = & \Phi(\xi_L,\xi_R) \derr^\Phi X(\xi_L,\xi_R) \Phi^\dagger(\xi_L,\xi_R), \\
\nonumber n^A_{jk}(\xi_L,\xi_R) & = & \Phi(\xi_L,\xi_R) \tilde A_{jk}(\xi_L,\xi_R) \Phi^\dagger(\xi_L,\xi_R), \\
\nonumber n^B_{jk}(\xi_L,\xi_R) & = & \Phi(\xi_L,\xi_R) \tilde B_{jk}(\xi_L,\xi_R) \Phi^\dagger(\xi_L,\xi_R),  \\
 n^C_{p}(\xi_L,\xi_R) & = & \Phi(\xi_L,\xi_R) \tilde C_{p} \Phi^\dagger(\xi_L,\xi_R). \label{movframeN}
\end{eqnarray}
satisfies the normalization conditions (\ref{normnorm}) and consequently the Gauss--Weingarten equations (\ref{gweqN}).


The Gauss--Codazzi--Ricci equations are the compatibility conditions for the Gauss--Weingarten equations (\ref{gweqN}) and are easily obtained comparing and equating the mixed derivatives, e.g. $\derr ( \derl {\cx}_L )$ and  
$\derr (\derl {\cx}_L)$, $\derl ( \derr n_j)$ and  $\derr ( \derl n_j)$ etc. They are the necessary and sufficient conditions for the local existence of the corresponding surface $\cf$. They are satisfied for any solution $P$ of the Euler--Lagrange equations (\ref{ELP}), provided that explicit forms of $Q^D_j,H_j$ are inserted (they can be found e.g. by the method developed below).

The second fundamental form and the mean curvature vector of the surface $\cf$ at the regular point $p$ can be expressed, according to \cite{Kob,Wil}, as
\bea
\nonumber {\bf II}  =  (\derl \derl X)^{\perp} {\rm d} \xi_L {\rm d} \xi_L + 2 (\derl \derr X)^{\perp} {\rm d} \xi_L {\rm d} \xi_R
+(\derr \derr X)^{\perp} {\rm d} \xi_R {\rm d} \xi_R, & & \\
\label{HN1}
{\bf H}  =  \frac{1}{\det G} \left( J_{R} (\derl \derl X)^{\perp}  - 2 G_{LR} (\derl \derr X)^{\perp} 
+ J_{L} (\derr \derr X)^{\perp} \right), & &
\eea
where $(\;)^{\perp}$ denotes the normal part of the vector.  
In our case 
\begin{eqnarray}
\nonumber 
(\derl \derl X)^{\perp} & = & [\derl \derl P,P] +
\frac{p_{L|R} \ p_{LL|R} }{ 4 - (p_{L|R})^2 } [\derl P,P] \\
\nonumber & - & \frac{ 2 p_{LL|R} }{ 4- (p_{L|R})^2 } [\derr P,P], \\
\nonumber (\derr \derr X)^{\perp} & = & - [\derr \derr P,P] + 
\frac{2 p_{RR|L} }{ 4- (p_{L|R})^2 } [\derl P,P] \\
\nonumber & - & \frac{ p_{L|R} \ p_{RR|L} }{ 4 - (p_{L|R})^2 } [\derr P,P], \\
(\derl \derr X)^{\perp} & = & \derl \derr X = [\derl P, \derr P].
\end{eqnarray}

\section{Reduction to the sine--Gordon equation}

In the special case of $n=m=1$
$$ G(1,1) \simeq \C P^1,$$
the formula (\ref{Gaussian}) simplifies to
$$ K=-4$$
and the mean curvature becomes
$$ H=2 i \frac{1+\tr(\derl P P \derr P)^2}{1-\tr(\derl P P \derr P)^2}$$

The fact that the $\C P^{1}$ sigma models lead to surfaces with constant negative Gaussian curvature suggests that there is an underlying connection between the $\C P^{1}$ sigma model and the sine--Gordon equation. Its explicit form can be most easily found from the first fundamental form written in the coordinates such that $J_L=1, \; J_R=1$, i.e. the equation (\ref{Iinconf})
$$ I = ({\rm d}\xi_L)^2 -  p_{L|R} {\rm d}\xi_L {\rm d}\xi_R +  ({\rm d} \xi_R)^2, $$
which can also be expressed as
$$ I = ({\rm d}\xi_L)^2 -  \left( \tr ( \derl P P \derr P ) + \tr ( \derr P P \derl P  ) \right) 
{\rm d}\xi_L {\rm d}\xi_R +  ({\rm d} \xi_R)^2. $$
The second fundamental form in such coordinates takes the form
$$ II = 2 ( [\derl P, \derr P] )  {\rm d}\xi_L {\rm d}\xi_R. $$
We can easily find a unit normal to the surface expressed as 
$$ n = -i (\1-2 P)$$
and the second fundamental form written as a scalar expression becomes
\bea
\nonumber II  & = & \frac{2i}{2} \tr \left( [\derl P,\derr P] (1-2 P) \right) {\rm d}\xi_L {\rm d}\xi_R \\
\nonumber & = & \frac{2}{i} \left( - \tr ( \derl P P \derr P )  + \tr ( \derr P P \derl P  ) \right) {\rm d}\xi_L {\rm d}\xi_R.
\eea
Altogether we see that the fundamental forms are given by 
\bea
\nonumber I & = & ({\rm d}\xi_L)^2 + 2 \cos(\phi)   {\rm d}\xi_L {\rm d}\xi_R +  ({\rm d} \xi_R)^2, \\
II  & = & 4 \sin(\phi) {\rm d}\xi_L {\rm d}\xi_R 
\eea
where
\be
{\rm e}^{i \phi} = - \tr ( \derl P P \derr P ).
\ee
From the differential geometry of surfaces it is immediately clear that $\phi$ satisfies the rescaled sine--Gordon equation
$$ \derl \derr \phi = 4 \sin \phi .$$
This can be also checked by a rather lengthy but straightforward explicit calculation. One should notice that the ``normalization'' of the solution 
$$J_L=1, \; J_R=1$$
is important -- otherwise a similar construction is possible, but a term proportional to $\sqrt{J_L J_R}$ appears on the right--hand side of the sine--Gordon equation. Luckily, in all the solutions of $\C P^{1}$ sigma model constructed below $J_L,J_R$ are constant so that the needed transformation of independent variables amounts to rescaling only.

To bring the sine--Gordon equation to its standard form one may of course rescale the independent variables, i.e.  define 
\be\label{indepvars}
\eta_L =  2 \xi_L, \qquad \eta_R  = 2 \xi_R.
\ee
Then 
\be\label{SGlc}
 \partial_{\eta_L} \partial_{\eta_R} \phi =  \sin \phi.
\ee
Defining 
$$ X= \eta_L + \eta_R, \qquad T=\eta_R- \eta_L $$
one gets the usual form of the sine--Gordon equation
\be\label{SG}
\partial_{TT} \phi-  \partial_{XX} \phi +  \sin \phi =0.
\ee
Finally, we use the Lorentz boost transformation
\be\label{Lorentz}
\tilde X= \frac{X-V T}{\sqrt{1-V^2}} ,\qquad \tilde T = \frac{T-V X}{\sqrt{1-V^2}}
\ee
with a suitable velocity $V$ in order to transform away the uniform motion of the wave and present only the significant properties of the solution. Such a transformation amounts to a rescaling
$$ \tilde \xi_L = \alpha \xi_L, \qquad \tilde \xi_R = \frac{1}{\alpha} \xi_R, \qquad \alpha = \sqrt{\frac{1+V}{1-V}}.$$

This relation between the $\C P^{1}$ sigma model and the sine--Gordon equation is not only of theoretical interest. It can be also used to construct nontrivial solutions of the sine--Gordon equation when solutions of $\C P^{1}$ are found. On the other hand, it can be also helpful in distinguishing qualitatively
different solutions of the $\C P^1$ sigma models, because solutions of the sine--Gordon equation are easier to visualize.

We should also mention that the relation is more precisely between the classes of equivalent solutions of the 
$\C P^1$ sigma model and of the sine-Gordon equation. Namely, in the $\C P^1$ case we had completely fixed the conformal transformations (\ref{conf}) by the requirement
$$J_L=1, \qquad J_R=1 .$$
This is in principle over--restrictive. For the transition to the sine--Gordon equation it would be sufficient to require that
$$J_L={\rm const.}>0, \qquad J_R=1/J_L.$$
However, the transformation between these choices of $J_L,J_R$ amounts to a Lorentz transformation on the sine--Gordon side, so that we are not losing anything provided we consider classes of Lorentz equivalent solutions of the sine--Gordon equation.

In the following we present three different kinds of solutions, namely
\begin{itemize}
\item Solutions of the $\C P^{1}$ sigma model obtained via symmetry reduction,
\item 1--soliton solutions of the $\C P^{1}$ sigma model obtained by the method of J. Harnad, Y. Saint-Aubin et al. \cite{HSA}.
\item 1--soliton solutions of the $\C P^{1}$ sigma model obtained by B. Piette \cite{Pie},
\end{itemize}

\subsection{Symmetry reduction of the Euler--Lagrange equation of the $\C P^1$ model}\label{symredsect}

For the purpose of the investigation of symmetries it is useful to rewrite the equation (\ref{ELP}) in terms
of the complex--valued function $w$, 
$$ P = \frac{1}{1+w \wc} \left( \ba{cc} w \wc & -\wc \\ -w  & 1 \ea \right).$$
The equation (\ref{ELP}) is then equivalently expressed in terms of $w$ as
\be\label{eqnmot}
\derl \derr w = 2 \frac{\wc \derl w \derr w}{1+w \wc},
\ee
and the quantities $J_L,J_R$ are
\be\label{CP1jljr}
 J_L = \frac{\derl w \derl \wc}{(1+w\wc)^2}, \; J_R = \frac{\derr w \derr \wc}{(1+w\wc)^2}.
\ee
The condition $J_L=1,J_R=1$ now becomes 
$$  \derl w \derl \wc = \derr w \derr \wc = (1+w\wc)^2$$
and the corresponding sine--Gordon solution $\phi$ is defined by 
$$ {\rm e}^{i \phi} = - \frac{\derl w}{\derr w}.$$

Writing (\ref{eqnmot}) in terms of real and imaginary parts
$$ w = u +i v.$$ 
we obtain
\begin{eqnarray}
\nonumber \derl \derr u  & = &
 \frac{2}{1+u^2+v^2} \left( u(\derl u \derr u - \derl v \derr v) + v(\derl v \derr u + \derl u \derr v) \right), \\
\nonumber \derl \derr v  & = &
 \frac{2}{1+u^2+v^2} \left( u(\derl v \derr u + \derl u \derr v) - v(\derl u \derr u - \derl v \derr v)  \right).
\label{eqnmotri}
\end{eqnarray}
The algebra of symmetry generators is infinite dimensional and can be expressed as the direct sum
\be\label{symalg}
{\cal G} = {\cal C}_{\xi_L} \oplus {\cal C}_{\xi_R} \oplus su(2),
\ee 
where ${\cal C}_{\xi_D},D=L,R$ denote infinite dimensional algebras of conformal transformations
$${\cal C}_{\xi_D}= \{ f_D(\xi_D) {\partial_{\xi_D}} | f_D \in {\cal C}^{\infty}(\R) \}$$
and $su(2)$ is generated by the following transformations involving only dependent coordinates
\begin{eqnarray}
\nonumber L_1 & = & u \partial_v - v \partial_u, \\
\nonumber L_2 & = & \frac{1}{2} (1+u^2-v^2) {\partial_u} + u v \partial_v, \\
 L_3 & = & - u v \partial_u +\frac{1}{2} (-1+u^2-v^2) \partial_v.
\end{eqnarray} 
For the construction of solutions invariant under some 1--parametric subgroup the conformal factors $f_D(\xi_D) $ can be absorbed into a suitable choice of independent variables, so that we are free to consider only one generator in each ${\cal C}_{\xi_D}$,
$$ \Xi_L = \partial_{\xi_L}, \; \Xi_R = \partial_{\xi_R}. $$
One finds that all solutions invariant under
$$ a \Xi_L + b \Xi_R $$
are singular, i.e. not usable for the construction of the associated surface.
Therefore we have to consider a solution invariant under a vector field which is a combination of $\Xi_L,\Xi_R,L_k$.
Since a vector from $su(2)$ itself cannot be used for symmetry reduction, because its orbits are not of codimension 1 in the space of independent variables, we can fix the $su(2)$ part of a general vector using the $SU(2)$ symmetry to be $L_1$ and consider only the following vector field 
\be\label{symredvec}
Y = L_1 +  a \Xi_L + b \Xi_R, \; a,b \in \R.
\ee
Using the method of characteristics one finds that a solution invariant under (\ref{symredvec}) must be of the form
\be\label{symredsoln}
 w = R(\chi) {\rm e}^{\frac{i}{a}(\xi_L-f(\chi))}, \; \; \chi=\xi_L-\frac{a}{b} \xi_R,
\ee
where $R,f: \R \rightarrow \R$. Substituting this form of $w$ into the Euler--Lagrange equation (\ref{eqnmot}) 
one finds two coupled ordinary differential
equations
\begin{eqnarray}
\label{ode1} R''-\frac{2 R}{1+R^2} {R'}^2 + \frac{R(1-R^2)}{1+R^2} (f'-{f'}^2) & = & 0,\\
\label{ode2} f''+ \frac{1-R^2}{R(1+R^2)} (2 R' f' -R') & = & 0,
\end{eqnarray}
where $R',f'$ etc. denote derivatives with respect to $\chi$. The system (\ref{ode1}),(\ref{ode2}) has a form similar to the one obtained by the symmetry reduction of the equations of the $\C P^{1}$ sigma model in (1+2)--dimensions in \cite{GWZ}. We shall now proceed in a way analogous to \cite{GWZ}. 

By integrating we rewrite (\ref{ode2}) in an equivalent form 
\be\label{ode2a}
f' = A \frac{(1+R^2)^2}{R^2}+\frac{1}{2},
\ee
where $A$ is a constant of integration. Substituting (\ref{ode2a}) into (\ref{ode1}) we get a single second order ODE
\be\label{ode1a}
R''-\frac{2 R}{1+R^2} {R'}^2 - A^2 \frac{(1-R^2)(1+R^2)^3}{R^3} + \frac{R(1-R^2)}{4(1+R^2)}  = 0.
\ee
Analyzing the singularity structure of the equation (\ref{ode1a}) we find that we can transform it into one of the standard Painlev\'e forms listed in \cite{Inc}. Performing the change of the dependent variable 
\be\label{Rsubs}
R(\chi) = \sqrt{-U(\chi)}
\ee
we find that the function $U$ obeys the Painlev\'e equation $\rm P_{XXXVIII}$ 
\be\label{XXXVIII}
U'' = \left( \frac{1}{2U} + \frac{1}{U-1} \right) {U'}^2+2 A^2 \frac{(1+U)(1-U)^3}{U}+\frac{U(1+U)}{2(U-1)}.
\ee
The order of (\ref{XXXVIII}) can be reduced by integration
\be\label{eqnU}
{U'}^2=-4 A^2U^4+4K U^3 + (8A^2-8K-1) U^2+4K U -4 A^2, \; \; K \in \C.
\ee
A considerable number of solutions of (\ref{eqnU}) exists \cite{GWZ}, but unfortunately most of them result either 
in $R$ being a complex function or in $f$ not being expressible in terms of elementary functions. Taking into account that  $R$ and $f$ are required to be real functions, we find the following solutions. (The list is not exhaustive, other solutions expressible in terms of elliptic integrals do exist. Due to their complexity they will be investigated in another paper.) The listed solutions have been already presented in proceedings \cite{Grusnosym}, the additional information presented here is the corresponding sine--Gordon reduction.

\subsubsection{The $\tanh$ solution}

As a first example of the construction of a surface let us consider
a special solution of (\ref{eqnU})
$$ U = - \tanh^2 \left( \frac{\chi-c}{4a} \right).$$
Consequently, we find from (\ref{Rsubs}), (\ref{ode2a}) that
\be\label{specsol}
R(\chi)=\tanh \left( \frac{\chi-c}{4a} \right), \; \; f(\chi)=\frac{\chi+d}{2},
\ee
$d \in R$ being a constant of integration . Finally, substituting (\ref{specsol}) 
into (\ref{symredsoln}) we find the solution of the Euler--Lagrange equation (\ref{eqnmot})
\be\label{example}
w =  \tanh \alpha \; {\rm e}^{i \beta},
\ee
where 
$$ \alpha = \frac{1}{4} \left( \frac{\xi_L}{a} - \frac{\xi_R}{b} - c \right), \; 
\beta = \frac{1}{2} \left( \frac{\xi_L}{a} + \frac{\xi_R}{b} - d \right) $$
and $a,b,c,d$ are real parameters.

In order to satisfy $J_L=1, J_R=1$ we rescale $\xi_L,\xi_R$, effectively putting
\be\label{abtanh} a= \pm \frac{1}{4}, \; b= \pm \frac{1}{4}.\ee
We choose the origin of our coordinates so that $c=d=0$.

\begin{figure}
\epsfxsize=5in
\begin{center}       
\leavevmode 
\epsffile{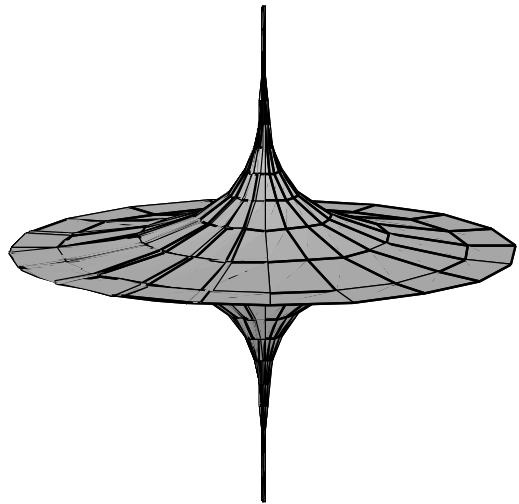}
\end{center}
\caption{Surface (\ref{tanhsurf}) associated with the $\tanh$ solution (\ref{example})} 
\label{tanhsurffig}
\end{figure}

Using the formula (\ref{ELPeqcons}) the corresponding surface is a Beltrami type pseudosphere immersed in $\R^3$ (see Fig. \ref{tanhsurffig}. Note that all figures presented in this paper were constructed using Maple 9 computer algebra system.)
 and can be written in the parametric form
\begin{eqnarray}
\nonumber X_1 & = & \frac{-\cos \beta}{2\cosh 2\alpha}+\frac{1}{2 \cosh 2}, \\
\nonumber X_2 & = & -\frac{\sin \beta}{2\cosh 2\alpha} ,\\
 X_3 & = & \frac{\tanh 2\alpha - \tanh 2 }{2}+1-\alpha. \label{tanhsurf}
\end{eqnarray}
The surface is shown in Figure \ref{tanhsurffig}.

The corresponding solution of the sine--Gordon equation (\ref{SG}) is 
\be
\phi = -4 \arctan \left( \epsilon_1 \tanh \frac{\eta_L+\epsilon_2 \eta_R}{2} \right)
\ee
where the signs $\epsilon_1,\epsilon_2$ depend on the choice of signs in (\ref{abtanh}). 
This represents a kink, going from $2 \pi$ to $0$ or vice versa.

\subsubsection{Exponential well solution}

\begin{figure}
\epsfxsize=5in
\begin{center}       
\leavevmode 
\epsffile{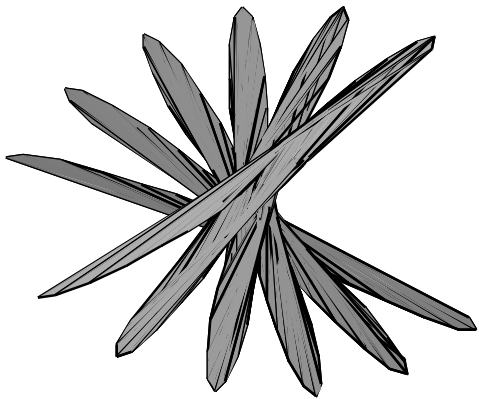}
\end{center}
\caption{The surface associated with the exponential well solution (\ref{expwell})} 
\label{expwellsurffig}
\end{figure}

As an example of an exponential solution of (\ref{ode1}),(\ref{ode2}) we select the following, so--called exponential well solution 
\begin{eqnarray}
\label{expwell}
 R ( \chi )  &=& \sqrt{\frac{(p-1) \cosh(g(\chi))+(p+1)}{(p-1) \cosh(g(\chi))-(p+1)} }, \\
\nonumber f(\chi) & = & 
\arctan \left( \frac{p+1}{2\sqrt{-p}} \tanh g(\chi) \right) +\frac{( p + 2 \sqrt{-p} -1) \chi  -2 \sqrt{-p} \chi_0 }{2(p-1)}+d 
\end{eqnarray}
where
$$ g(\chi)=\frac{(p+1)(\chi-\chi_0)}{2(p-1)}, \; p<-1. $$
The solution of the Euler--Lagrange equation (\ref{eqnmot}) is expressed using the formula (\ref{symredsoln})
$$ w = R(\chi) {\rm e}^{i(\xi_L/a-f(\chi))}, \; \; \chi=\frac{\xi_L}{a}-\frac{\xi_R}{b}. $$

The condition $J_L=1,J_R=1$ requires a rescaling of $\xi_L,\xi_R$ that fixes $a,b$
$$ a = \frac{p-2 \sqrt{-p}-1}{4 (p-1)}, b= \frac{p+2 \sqrt{-p}-1}{4 (p-1)}. $$

A picture of the surface is given in Fig. \ref{expwellsurffig} for the values of parameters
$$ p=-\frac{3}{2},\;\chi_0=0, \; d=0, \; \xi_L,\xi_R \in (-40 \ldots 40).$$
This example represents an immersed Willmore surface that is not a stereographic projection of a compact immersed minimal surface in $S^3$. 

The sine--Gordon solution is again a kink. The value of the parameter $p$ determines its steepness.

\subsubsection{Elliptic solution}

\begin{figure}
\epsfxsize=5in
\begin{center}       
\leavevmode 
\epsffile{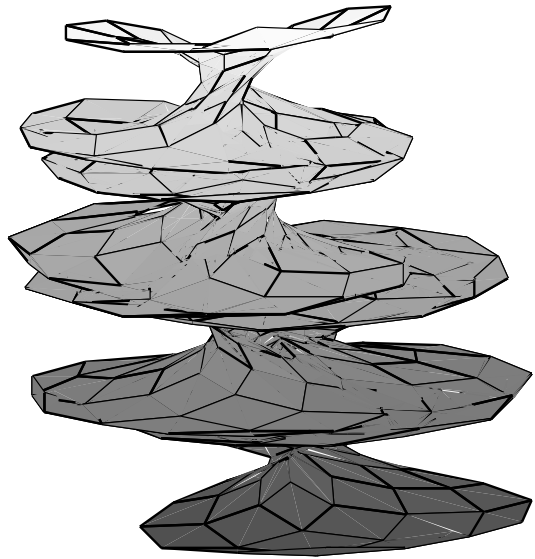}
\end{center}
\caption{The surface associated with the elliptic solution (\ref{ellsol}) is a type of a pseudospherical 
helicoid} 
\label{ellsolfig}
\end{figure}

\begin{figure}
\epsfxsize=3in
\begin{center}       
\leavevmode 
\epsffile{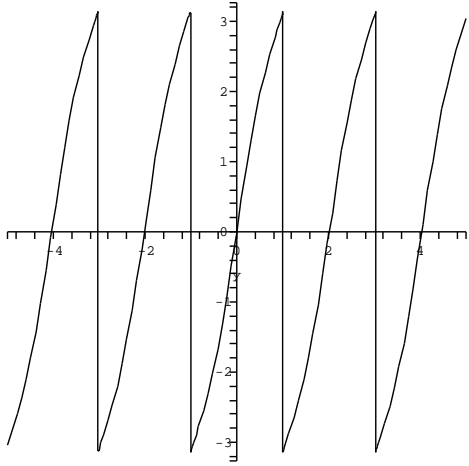}
\end{center}
\caption{The associated static solution of the sine--Gordon equation as a function of the space coordinate $X$} 
\label{ellsolsGfig}
\end{figure}

There exists also a class of solutions of (\ref{ode1}),(\ref{ode2}) which can be written in terms of elliptic functions. We select for the construction of a surface one of these solutions which is written in terms of Jacobi sn function 
\be
R( \chi ) =                                            
\sqrt{-p} \ {\rm sn} \left( \sqrt{K q} (\chi_0-\chi), \sqrt{ \frac{ p }{ q } } \right), \; \; f(\chi)=\frac{\chi+d}{2}
\ee
where
$$ p=\frac{1+8 K-\sqrt{1+16 K}}{8K}, \; q= \frac{1+8 K+\sqrt{1+16 K}}{8K}$$
and in order for $R(\chi)$ to be real
$$K \in (-\frac{1}{16},0).$$
The solution of the Euler--Lagrange equation (\ref{eqnmot}) is therefore 
\be\label{ellsol}
w(\xi_L,\xi_R)=\sqrt{-p} \ {\rm sn} \left( \sqrt{K q} (\xi_0-\frac{\xi_L}{a}+\frac{\xi_R}{b}), \sqrt{ \frac{ p }{ q } } \right) 
{\rm e}^{\frac{i}{2} \left( \xi_L/a+\xi_R/b - d \right)}.
\ee
The $J_L=1,J_R=1$ condition fixes $a,b$
$$ a = \frac{p-2 \sqrt{-p}-1}{4 (p-1)}, b= \frac{p+2 \sqrt{-p}-1}{4 (p-1)}. $$
We present a picture of the associated surface in Fig. \ref{ellsolfig} for the parameters
$$ K=-\frac{1}{20}, \; \xi_L,\xi_R \in (-10,\ldots,10). $$

The sine--Gordon solution $\phi$ is then defined via
$$ \phi= -I \ln \left( -\frac{Q {\rm cn}(Q X,\mu) {\rm dn}(Q Y,\mu)-2i {\rm sn}(Q Y,\mu)}{ Q {\rm cn}(QY,\mu) {\rm dn}(Q Y,\mu)+2i {\rm sn}(Q Y,\mu) } \right) $$
where
$$ Y=\frac{1}{8\sqrt{-K}} X, \; Q=\sqrt{2+16 K+2 \sqrt{1+16K}}, \; \mu=\sqrt{\frac{1+8K-\sqrt{1+16K}}{1+8K+\sqrt{1+16K}} } .$$
The solution is a static generalization of a kink, the value of the parameter $K$ determines the steepness of the solution, discontinuities arise from $\ln$ involved in the construction of $\phi$ (see Fig. \ref{ellsolsGfig}).

\subsection{1--soliton solutions obtained via the method of J. Harnad, Y. Saint-Aubin et al.}\label{HSAsect}

When constructing 1--soliton solutions of the $\C P^{1}$ sigma model using the method described in \cite{HSA}
one starts from a trivial, vacuum solution which is represented by
\be
 g = \left( \ba{cc} \cos t &  \sin t \\ -\sin t & \cos t \ea \right), \; t=\frac{\xi_L+\xi_R}{2} 
\ee
where the corresponding projector is written as
\be\label{jsolP}
 P = \frac{1}{2} \left( 1- g .{\rm diag}(1,-1) \right).
\ee
Then one finds a matrix $\psi(\lambda)$ such that
$$ \derl \psi(\lambda) = \frac{1}{1+\lambda} (\derl g g^{-1}) \psi, \; \derr \psi = \frac{1}{1-\lambda} (\derr g g^{-1}) \psi,$$
(note that $\partial_{L,R} g g^{-1}$ are  constant for the given special $g$).
A solution to this equation is given by
$$ \psi(\lambda) = \left( \ba{cc} \psi_{11} & \psi_{12} \\ \psi_{21} & \psi_{22} \ea \right) .$$ 
where
\bea 
\nonumber \psi_{11} & = & \cos \frac{\xi_R}{2(\lambda-1)} \cos \frac{\xi_L}{2(1+\lambda)}
+\sin \frac{\xi_L}{2(1+\lambda)} \sin \frac{\xi_R}{2(\lambda-1}, \\
\nonumber  \psi_{12} & = &  -\sin \frac{\xi_R}{2(\lambda-1)} \cos \frac{\xi_L}{2(1+\lambda)}
+\sin \frac{\xi_L}{2(1+\lambda)} \cos \frac{\xi_R}{2(\lambda-1)}, \\
\nonumber  \psi_{21} & = & \sin \frac{\xi_R}{2(\lambda-1)} \cos \frac{\xi_L}{2(1+\lambda)}
- \sin \frac{\xi_L}{2(1+\lambda)} \cos \frac{\xi_R}{2(\lambda-1)}, \\
\nonumber  \psi_{12} & = & \cos \frac{\xi_R}{2(\lambda-1)} \cos \frac{\xi_L}{2(1+\lambda)} 
+ \sin \frac{\xi_L}{2(1+\lambda)} \sin \frac{\xi_R}{2(\lambda-1)}. 
\eea
Further defining
$$ M = \psi(\bar\lambda) \left( \ba{c} \alpha \\ \beta \ea \right),$$
$$ R = M ( M^\dg M )^{-1}  M^\dg, $$
$$ U= \1 + ( \bar \lambda-\lambda)/\lambda  R$$
$$ \tilde g = U g $$
one obtains a new $g$ which is a solution of the same equation as the original $g$, namely
$$ \derr (\derl g g^{-1}) + \derl (\derr g g^{-1}) =0 .$$
If we further impose the conditions 
$$ |\lambda|=1, \; |\alpha|=|\beta| $$
the new solution can again be written in terms of a projector like in (\ref{jsolP}), i.e. $\tilde g$ defines another 
solution of the $\C P^{1}$ sigma model.

We have used the method given above and constructed the corresponding 1--soliton solutions, using the computer algebra system Maple. Unfortunately, the expressions which we obtained are extremely complicated and even the solutions ($g$ or $P$) are impossible to display. Since they can be easily recalculated using the process explained above, we only present one figure of the obtained surface (integrated numerically, of course).

\begin{figure}
\epsfxsize=5in
\begin{center}       
\leavevmode 
\epsffile{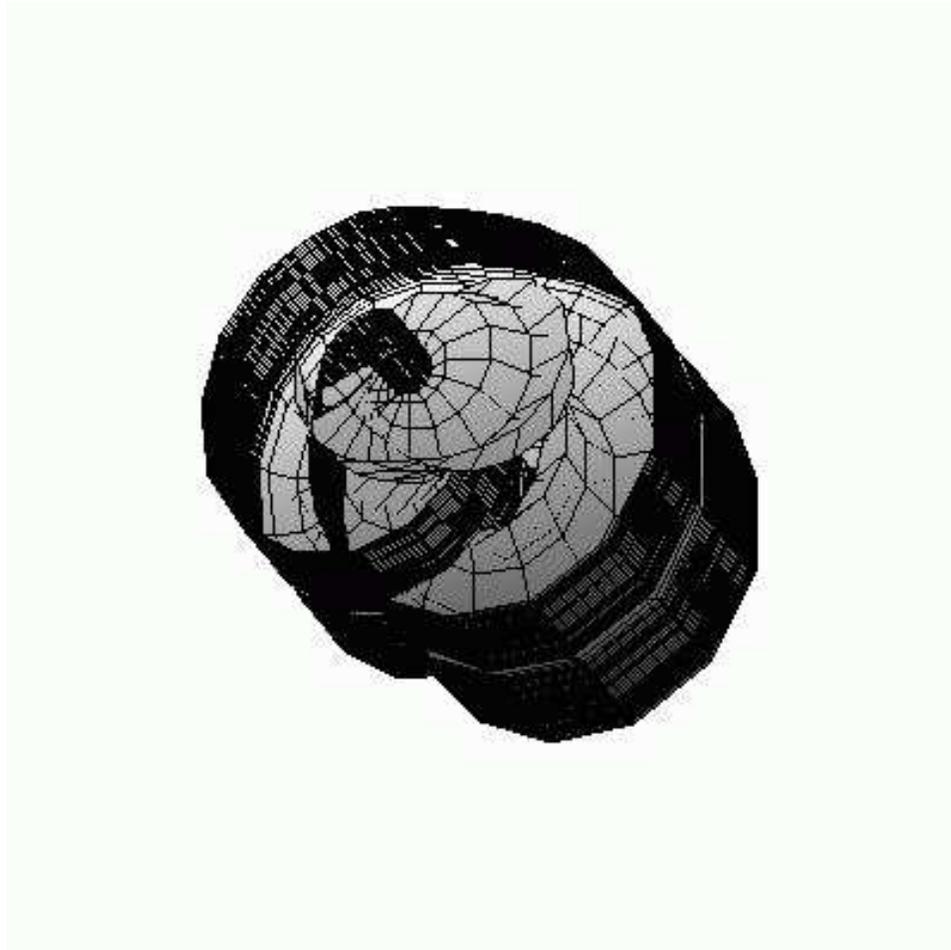}
\end{center}
\caption{The surface associated with 1--soliton solution for 
$\lambda=\frac{1}{\sqrt{2}}+\frac{1}{\sqrt{2}} i, \; \alpha=-\sqrt{2}, \beta=1+i $ appears to be a Dini type surface} 
\end{figure}

The corresponding solution of the sine--Gordon equation appears to be a kink. Unfortunately, the complicated form of the projector $P$ obtained in this way leads to serious difficulties in even a numerical computation of the corresponding solution of the sine-Gordon equation. Therefore, we were able to calculate graphs of only a few time slices (and even those took several days of computation in Maple) and from their comparison we extrapolated the behaviour of the solution. In the case described in the next section the analysis was finally also preformed numerically but we were able to compute the solution for a sufficiently large number of points so that we can be confident of the validity of our results in that case. 

\subsection{1--soliton solutions constructed by B. Piette and the corresponding surfaces and solutions of the sine--Gordon equation}

In \cite{Pie} B. Piette constructs a 1--parameter family of solutions of the $\C P^{1}$ sigma model 
which can be written in the following form
\be\label{psol}
P = \frac{1}{2} \left( \ba{cc} 1 - g_{11} & -g_{12} \\ -\bar g_{12} & 1 + g_{11}  \ea \right) 
\ee
where
\bea  
\nn \sigma  & = & 2 (\xi_L+\xi_R), \; \; \tau  =  2 \left( \frac{\xi_L}{1+\bar \lambda}+\frac{\xi_R}{1-\bar \lambda} \right), \\
\nn u & = & 2 \ \Re \ \tau,  \; v=2 \ \Im  \ \tau, \\
\nn \Lambda^2 & = & \bar \lambda \lambda \left( \frac{4 \cos^2 (u-\sigma)}{(1-\bar\lambda \lambda)^2}+\frac{\cosh^2(v)}{(\Im \lambda)^2} \right), \\
\nn  g_{11} & = & \cos \sigma -\frac{1}{\Lambda^2} 
\left( \frac{\Re \lambda \sin \sigma \sinh(2v)+\Im \lambda \cos \sigma \cosh(2v)}{\Im \lambda} \right. \\
\nn & + & \left. \frac{\cos(2 u-3\sigma)-\bar\lambda \lambda \cos(2u-\sigma)}{\bar\lambda\lambda-1} \right), \\
\nn g_{12}  & = & -\sin \sigma - \frac{1}{\lambda^2} \left[\frac{\Re \lambda\cos \sigma \sinh(2v)-\Im \lambda \sin \sigma \cosh(2v)}{\Im \lambda} \right. \\
\nn  & + & \frac{\sin(2u-3\sigma)+\bar\lambda\lambda\sin(2u-\sigma)}{\bar\lambda\lambda-1} \\
\nn & + & \left. 2i \left( \frac{\Re \lambda}{\Im \lambda} \sin(u-\sigma)\cosh v-\frac{\bar\lambda\lambda+1}{\bar\lambda\lambda-1} \cos(u-\sigma)\sinh v \right)
 \right].
 \eea
The complex parameter $\lambda$ is supposed to satisfy $|\lambda| \neq 1$. Compared to \cite{Pie} 
the coordinates were rescaled so that we have $J_L=J_R=1$.

In \cite{Pie} the behaviour of such solitons is briefly discussed. We shall present the corresponding surfaces and associated solutions of the sine--Gordon equation. Because the resulting formulae for tangent vectors etc. become rather complicated (and the Weierstrass representation of surfaces can be integrated only numerically), we resort to numerical calculations and present images of surfaces for a few values of $\lambda$. Rather surprisingly, the value of the parameter plays an essential role in the shape of the surface, for different $\lambda$s the surfaces look distinctly different. 

\begin{figure}
\epsfxsize=5in
\begin{center}       
\leavevmode 
\epsffile{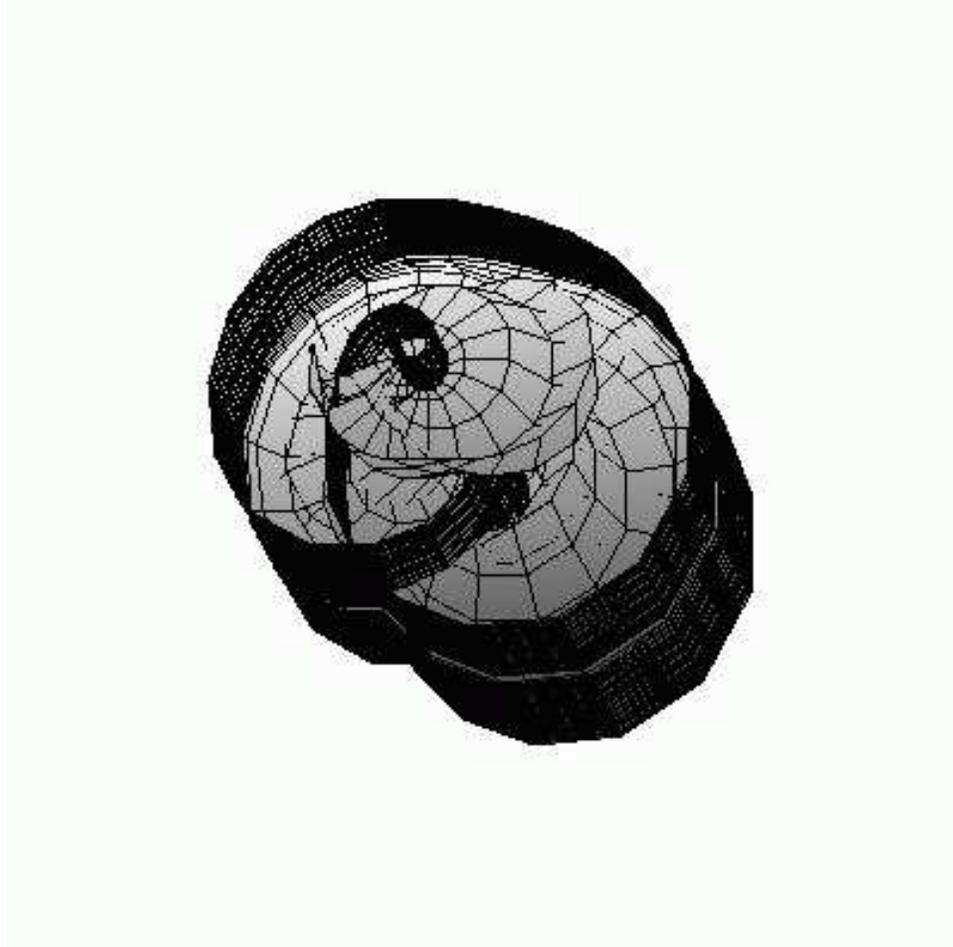}
\end{center}
\caption{The surface associated with $\lambda= 1.1+1.1i$. This is the choice of $\lambda$ originally considered in \cite{Pie}.} 
\end{figure}

\begin{figure}
\epsfxsize=3in
\begin{center}       
\leavevmode 
\epsffile{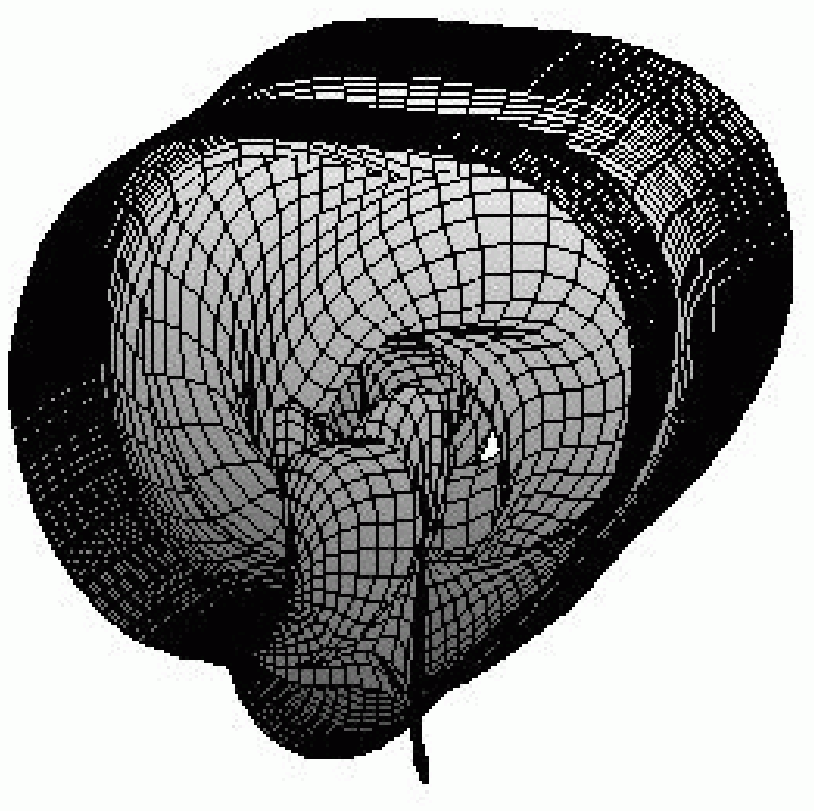}
\end{center}
\caption{The surface associated with $\lambda= 1+2i$} 
\end{figure}

\begin{figure}
\epsfxsize=5in
\begin{center}       
\leavevmode 
\epsffile{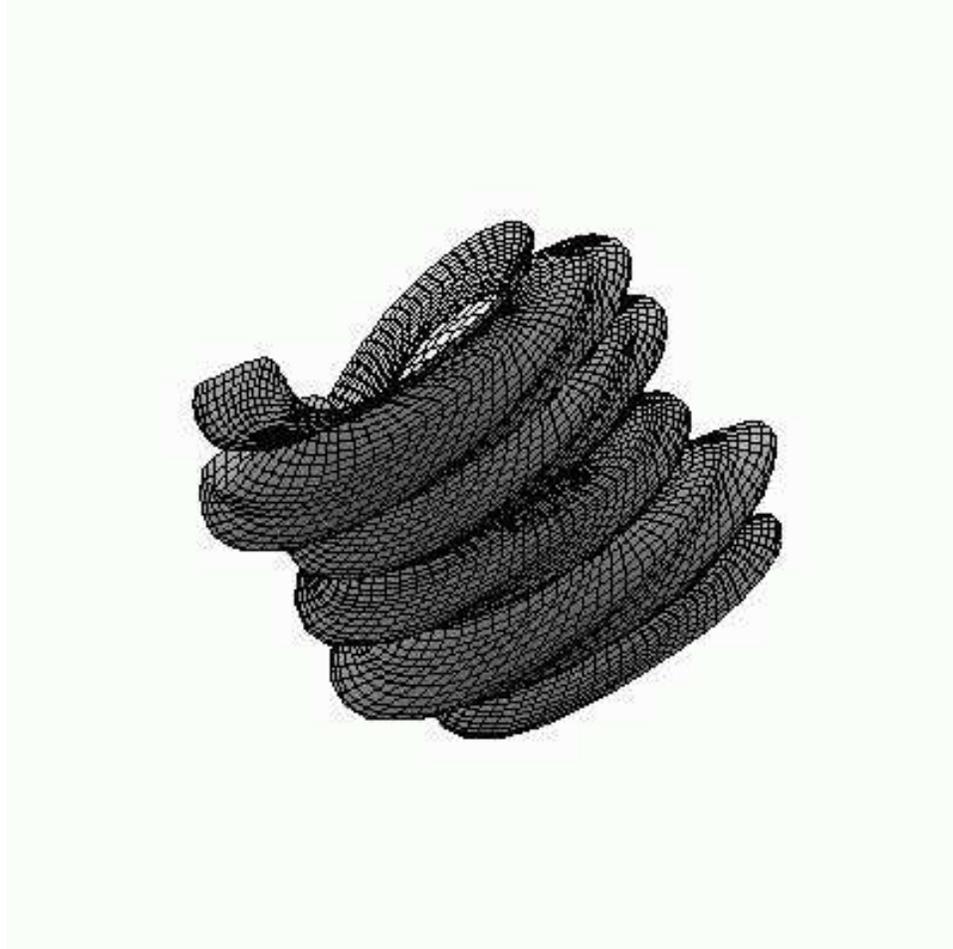}
\end{center}
\caption{The surface associated with $\lambda= \pi+\frac{1}{2}i$} 
\end{figure}

\begin{figure}
\epsfxsize=5in
\begin{center}       
\leavevmode 
\epsffile{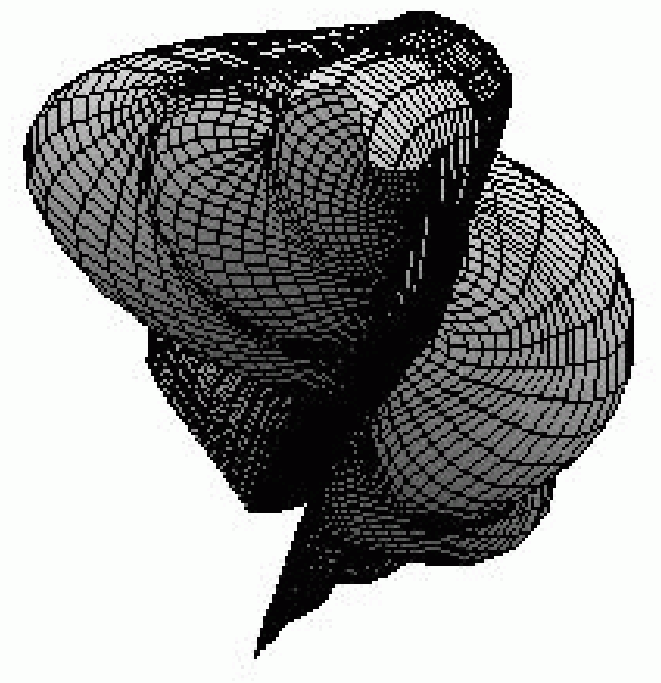}
\end{center}
\caption{The surface associated with $\lambda= -2+2i$} 
\end{figure}

\begin{figure}
\epsfxsize=3in
\begin{center}       
\leavevmode 
\epsffile{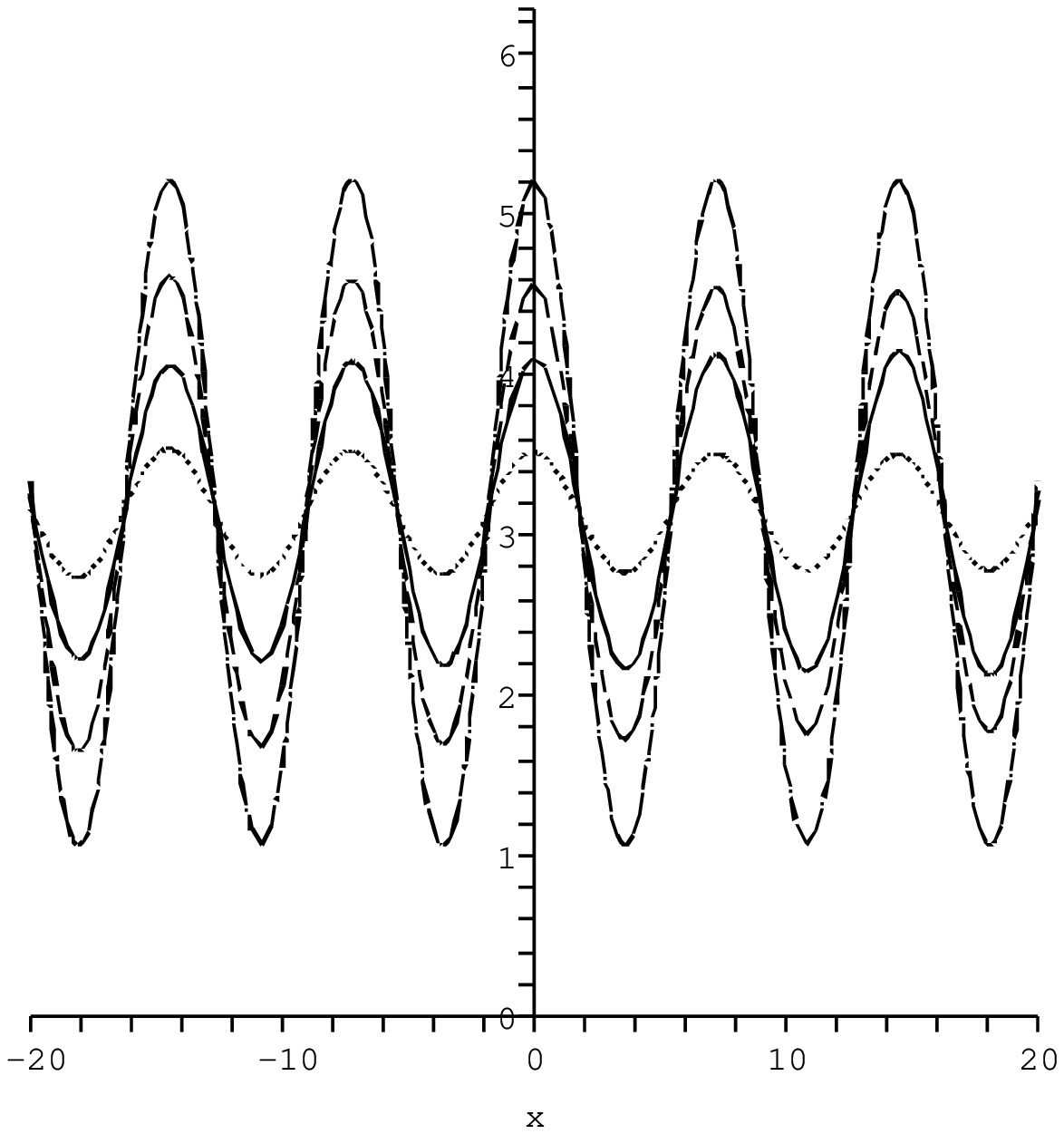}
\end{center}
\caption{$\lambda= -2+2i.$ The associated solution of the sine--Gordon equation -- time slices} 
\label{plm22sGts}
\end{figure}

The corresponding solutions of the sine--Gordon equation are periodic, stationary (after a properly chosen Lorentz boost). They can be constructed analytically using any computer algebra systems, but unfortunately the resulting formulae are too complicated to analyze. Therefore, we rely on the numerical computations of their graphs.
The graphs of the solutions look rather similar for all values of $\lambda$ we have investigated; only the width and height of $\phi$ appear to depend on $\lambda$. Therefore, we present only one example, namely for $\lambda= -2+2i$ in Fig. \ref{plm22sGts}. The solution is shown as a function of $X$ in $T=-3,0,2,5$ (solid, dotted, dashed, dotted--dashed line, respectively). 

Note that the sine--Gordon solutions in this case and in the subsection \ref{HSAsect} are different, and consequently the solutions of the $\C P^{1}$ sigma model and surfaces in these cases are necessarily also different.

\subsection{Linear problem for the $\C P^{1}$ model}
Finally we should mention that there is a hope that other solutions may be found using the Lax pair  
representation of the $\C P^{1}$ sigma model. We know two rather different Lax pairs. The first one 
arises in the construction of the moving frame. The group element $\Phi$ bringing the tangent vectors to the form
(\ref{Phireq}) satisfies
$$ \derl \Phi = \tilde U \Phi, \quad \derr \Phi = \tilde V \Phi $$
where $\tilde U, \tilde V\in su(2)$ are of the form
\begin{eqnarray}
\tilde U & = & \frac{1}{1+w \wc} \left(
\begin{array}{cc}
\frac{1}{2} (w \derl \wc - \wc \derl w)  & -\derl \wc  \\
\derl w  &  \frac{1}{2} (\wc \derl w - w \derl \wc)  \end{array} \right),   \nonumber \\
\tilde V & = & \frac{1}{1+w \wc} \left(
\begin{array}{cc}
\frac{1}{2} (w \derr \wc - \wc \derr w)  & -\derr \wc  \\
\derr w  &  \frac{1}{2} (\wc \derr w - w \derr \wc)  \end{array} \right). \nonumber 
\end{eqnarray}
Consistency requires that
\be
\derr \tilde U - \derl \tilde V + [\tilde U,\tilde V] =0
\ee
Such $\tilde U,\tilde V$ can be generalized to $sl(2,\C)$ matrices involving a spectral parameter $\lambda$
\begin{eqnarray}
 \tilde U_\lambda  & = & \frac{1}{1+w \wc} \left(
\begin{array}{cc}
\frac{1}{2} (w \derl \wc - \wc \derl w)  & -\lambda \derl \wc  \\
\derl w  &  \frac{1}{2} (\wc \derl w - w \derl \wc)  \end{array} \right),   \nonumber \\
\tilde V_\lambda & = & \frac{1}{1+w \wc} \left(
\begin{array}{cc}
\frac{1}{2} (w \derr \wc - \wc \derr w)  & -\derr \wc  \\
\frac{1}{\lambda } \derr w  &  \frac{1}{2} (\wc \derr w - w \derr \wc)  \end{array} \right).  
\label{laxpair}
\end{eqnarray}
and satisfying 
\be\label{lincons}
\derr \tilde U_\lambda - \derl \tilde V_\lambda + [\tilde U_\lambda,\tilde V_\lambda] =0
\ee
for any solution $w$ of the Euler--Lagrange equation (\ref{eqnmot}) (note that (\ref{lincons}) is equivalent to 
$$ \frac{ (\lambda-1) \left( \derl \derr w (1+ w \wc) - 2 \wc \derl w \derr w \right)}{\lambda (1+w \wc)^2 }=0$$
and its complex conjugate.)

Consequently, the associated linear problem (Lax pair) reads 
\be\label{asslinsys}
 \derl \psi = \tilde U_\lambda \psi, \qquad  \derr \psi = \tilde V_\lambda \psi.
\ee
 
Another Lax pair for the equation (\ref{eqnmot}) can be  written in the form \cite{ZaMi}
\be\label{zamich}
 \derl \Psi = \frac{2}{1+\lambda} M_L \Psi, \; \derr \Psi = \frac{2}{1-\lambda} M_R \Psi.
\ee
It is interesting to note that the compatibility conditions for both linear systems (\ref{asslinsys})
and (\ref{zamich}) gives us exactly the equation (\ref{eqnmot}) but the structure of the Lax operators is different.
In both cases they involve  first order derivatives of $w,\wc$, but (\ref{zamich}) involves the spectral parameter $\lambda$ in an overall factor whereas (\ref{laxpair}) has $\lambda$ only in some of the matrix entries. This means that for some purposes one Lax pair may be more suitable that the other one.

\section{Final remarks}\label{concl}
We have presented a rather straightforward procedure for the construction of surfaces from solutions of the Grassmannian sigma models, and derived formulae describing their geometric properties like Gaussian curvature, structural equations etc. Compared to our previous works \cite{Grusno,Grusnosym,GrusnoGr}, the main improvement lies in the simpler formulation using projectors, which are in $1-1$ correspondence with elements of $G(m,n)$. This allows to  avoid from the beginning the superfluous gauge degrees of freedom involved in the use of representatives of elements of $G(m,n)$ as equivalence classes in $su(N)$, as was employed in \cite{GrusnoGr}, and consequently to get the better understanding of the procedure and compact expressions for the important objects like fundamental forms, curvatures etc. 

The approach presented here is from the beginning constrained by the specific choice of the field theory on two--dimensional spacetime, i.e. Grassmannian sigma models. A question naturally arises whether it can be generalized to other field theories, possibly on different spacetimes. For example, the $\C P^{N}$ sigma models in (2+1) dimensions may in principle lead to the construction of 3--dimensional submanifolds of specific properties immersed in $su(N+1)$. Similarly, it might be of interest in applications, especially in physics, to investigate whether a similar approach can be applied also in the case of non--Abelian field theories.

Other aspects of the method worth investigating are whether the method may under some condition lead to compact surfaces, whether the surfaces are stable under perturbation of the solution of the underlying sigma model etc.

In the $\C P^{1}$ case the most interesting open question is the following one. We have showed that any ``regular'', i.e. such that $\tr (\derl P.\derl P), \tr (\derr P.\derr P) \neq 0$, solution of the $\C P^{1}$ gives (locally) rise to a solution of the sine--Gordon equation. However, we don't know whether for every solution of the sine--Gordon equation there exists a corresponding solution of the $\C P^{1}$ sigma model or whether any useful criterion can be found determining which solutions of the sine--Gordon equation may be obtained in this way. 

These and other question are problems for future research.

\subsection*{Acknowledgements}

This work was supported in part by research grants from NSERC of Canada. Libor \v Snobl acknowledges a postdoctoral fellowship awarded by the Laboratory of Mathematical Physics of the CRM, Universit\'e de Montr\'eal. 
We have benefited from helpful discussions with Pavel Winternitz on the subject of this work.

\end{document}